\documentclass[12pt,oneside,reqno]{amsart}
\usepackage{mathrsfs}
\usepackage{graphics}
\usepackage{amssymb}
\usepackage{enumerate}
\newcommand{\dif}{\mathrm{d}}

\newcommand{\be}{\begin{eqnarray}}
\newcommand{\ee}{\end{eqnarray}}
\newcommand{\ce}{\begin{eqnarray*}}
\newcommand{\de}{\end{eqnarray*}}
\newtheorem{theorem}{Theorem}[section]
\newtheorem{lemma}[theorem]{Lemma}
\newtheorem{remark}[theorem]{Remark}
\newtheorem{definition}[theorem]{Definition}
\newtheorem{proposition}[theorem]{Proposition}
\newtheorem{Example}[theorem]{Example}
\newtheorem{corollary}[theorem]{Corollary}
\newtheorem{condition}[theorem]{Condition}
\def\e{\varepsilon}
\def\s{\sigma}
\def\t{\theta}
\def\a{\alpha}

\def\b{\beta}

\def\g{\gamma}

\def\[{{\Big[}}
\def\]{{\Big]}}
\def\<{{\langle}}
\def\>{{\rangle}}
\def\({{\Big(}}
\def\){{\Big)}}

\def\no{\nonumber}
\def\bt{\begin{theorem}}
\def\et{\end{theorem}}
\def\bl{\begin{lemma}}
\def\el{\end{lemma}}
\def\br{\begin{remark}}
\def\er{\end{remark}}
\def\bx{\begin{Example}}
\def\ex{\end{Example}}
\def\bd{\begin{definition}}
\def\ed{\end{definition}}
\def\bp{\begin{proposition}}
\def\ep{\end{proposition}}
\def\bc{\begin{corollary}}
\def\ec{\end{corollary}}
\def\bco{\begin{condition}}
\def\eco{\end{condition}}

\def\mE{{\mathbb E}}

\def\mN{{\mathbb N}}

\def\mP{{\mathbb P}}

\def\mR{{\mathbb R}}

\def\sF{{\mathscr F}}

\def\geq{\geqslant}
\def\leq{\leqslant}
\def\epsilon{\varepsilon}
\begin{document}

\allowdisplaybreaks

\title{The central limit theorem for stochastic Volterra equations with singular kernels}

\author{Huijie Qiao}

\dedicatory{School of Mathematics,
Southeast University,\\ Nanjing, Jiangsu 211189 China\\
hjqiaogean@seu.edu.cn}

\thanks{{\it AMS Subject Classification(2020):} 60H20; 60F05}

\thanks{{\it Keywords:} Stochastic Volterra equations with singular kernels; the central limit theorem; fractional Brownian motions}

\thanks{This work is supported by NSF of China (No. 12071071)}

\subjclass{}

\date{}

\begin{abstract}
This work concerns stochastic Volterra equations with singular kernels. Under the suitable conditions, we prove the central limit theorem for them. Moreover, we apply our result to stochastic Volterra equations with the kernels of fractional Brownian motions with the Hurst parameter $H\in(0, 1)$.
\end{abstract}

\maketitle \rm

\section{Introduction}
Fix $T>0$ and consider the type of stochastic Volterra equations 
\be
X_t(x)=x+\int_0^tK_1(t,s)b(s,X_s(x))\dif s+\int_0^tK_2(t,s)\s(s,X_s(x))\dif B_s, \quad t\in[0,T],
\label{meq0}
\ee
where $x\in\mR^d$, $K_i(t,s), i=1,2$  are two positive functions on $[0,T]\times[0,T]$ which may be singular, $(B_t)_{t\in[0,T]}$ is a $m$-dimensional Brownian motion defined on a filtered probability space $(\Omega,\sF, (\sF_t)_{t\in[0,T]}, \mP)$, and the coefficients $b: \mR_+\times\mR^d\mapsto\mR^d$, $\sigma: \mR_+\times\mR^d\mapsto\mR^d\times\mR^m$ are all Borel measurable. This kind of stochastic Volterra equations  appear in many fields, such as nonlinear filtering \cite{cd0}, fluid turbulence \cite{ch}, turbulence modelling in atmospheric winds \cite{chpp} and mathematical finance \cite{cp}. And there have been many related results (c.f. \cite{ajlp, g, jp, lwyz, nr, qw, wang, zhang}).

In this paper, we are devoted to observing the asymptotic behavior of small perturbation for Eq.(\ref{meq0}). Concretely speaking, consider the following stochastic Volterra equation
\be
X_t^\e(x)=x+\int_0^tK_1(t,s)b(s,X^\e_s(x))\dif s+\sqrt{\e}\int_0^tK_2(t,s)\s(s,X^\e_s(x))\dif B_s, \quad t\in[0,T].
\label{meq1}
\ee
We will study the asymptotic behavior of the amount
\be\label{ab}
\frac{X^{\epsilon}_{\cdot}- X^{0}_\cdot}{\sqrt{\epsilon}} ~\mbox{in}~ C([0,T]\times\mR^d, \mR^d)
\ee
as $\epsilon\to 0$, where $X^{0}$ solves the following Volterra equation:
\be
X^0_t(x)=x+\int_0^tK_1(t,s)b(s,X_s^{0}(x))\dif s.
\label{volequ}
\ee

For (\ref{ab}), if $\sqrt{\epsilon}$ is replaced by the constant $1$, the result for this class of asymptotic behaviors is called the large deviation principle (LDP for short). In \cite{nr}, Nualart and Rovira first studied the LDP for stochastic Volterra equations with regular kernels. On one hand, Lakhel \cite{la} improved this result to the Besov-Orlicz space. On the other hand, Zhang \cite{zhang} extended this result to stochastic Volterra equations with singular kernels. Later, for special Volterra type rough volatility models, Cellupica and Pacchiarotti \cite{cp} and Gulisashvili \cite{g} both established corresponding LDPs. Very recently, for general stochastic Volterra systems, Jacquier and Pannier \cite{jp} showed the LDP under weak conditions.

In (\ref{ab}), if $\sqrt{\epsilon}$ is replaced by the function $a(\e)$ which satisfies
    \ce
	a(\epsilon)\rightarrow0,\quad \frac{\epsilon}{a^{2}(\epsilon)}\rightarrow0\quad as\quad \epsilon\rightarrow0,
	\de
we call the result for this asymptotic behavior the moderate deviation principle (MDP for short). In \cite{lwyz}, Li et al. proved the MDP for stochastic Volterra equations with regular kernels in a $1$-dimensional space. In \cite{jp}, Jacquier and Pannier generalized this result to stochastic Volterra systems with singular kernels.

In the case of (\ref{ab}), the result for this class of asymptotic behaviors is called the central limit theorem (CLT for short). The CLT describes that $\frac{X^{\epsilon}- X^{0}}{\sqrt{\epsilon}}$ converges to a stochastic process in a certain sense as $\epsilon\rightarrow0$. That is, from the CLT, we know the concrete limit. However, up to now, there is yet no related result about stochastic Volterra equations. In this paper, we prove that under some suitable conditions, $\frac{X^{\epsilon}- X^{0}}{\sqrt{\epsilon}}$ converges to the solution of a stochastic Volterra equation in the $L^p$ sense as $\epsilon\rightarrow0$.

Notice that the solution of Eq.(\ref{meq1}) is in general not a semimartingale nor a Markov process, preventing the usage of It\^o calculus or Feynman-Kac type formulas. Therefore, in this paper, we apply a lot of techniques to obtain some estimates.

Finally, we formulate our motivation of this paper. Note that by integral transformation of stochastic integrals with respect to fractional Brownian motions, one can link them to the stochastic integrals with respect to standard  Brownian motions involving singular kernels, thus stochastic differential equations driven by fractional Brownian motions can be treated as stochastic Volterra equations involving singular kernels. This motives us in this paper to study the CLT for stochastic Volterra equations with singular kernels, with stochastic differential equations driven by fractional Brownian motions as an important special subclass. 

The paper is organized as follows. In the next section, we state the main result. In Section \ref{cltthproof}, the main theorem is proved in details. Finally, we show how the result apply to stochastic differential equations driven by fractional Brownian motions.

In the following $C$ with or without indices will denote different positive constants whose values may change from line to line.

\section{The central limit theorem for stochastic Volterra equations with singular kernels}\label{clt}

In this section, we study the CLT for stochastic Volterra equations with singular kernels.

First of all, for Eq.(\ref{meq1}) and Eq.(\ref{volequ}), we assume:
\begin{enumerate}[(${\bf H}_{K}^1$)]
\item There is $\b>1$ such that
\ce
\sup\limits_{t\in[0,T]}\int_0^t\left[K_1(t,s)^\b+K_2(t,s)^{2\b}\right]\dif s\leq C.
\de
\end{enumerate}
\begin{enumerate}[(${\bf H}_{K}^2$)]
\item There is $\g>0$ such that for any $t, t'\in[0,T]$
\ce
\int_0^{t\land t'}\left[|K_1(t',s)-K_1(t,s)|+|K_2(t',s)-K_2(t,s)|^2\right]\dif s\leq C|t'-t|^\g.
\de 
\end{enumerate}
\begin{enumerate}[(${\bf H}_{b,\s}^{1}$)]
	\item $b,\sigma$ satisfy for any $t\in[0,T]$ and $x, x_1, x_2\in\mR^d$
	\ce
	&&\|\nabla b(t,x)\|\leq L_1, \quad |b(t,x)|\leq L_1(1+|x|), \\
    &&\|\sigma(t,x_1)-\sigma(t,x_2)\|\leq L_1|x_1-x_2|,\quad \|\s(t,x)\|\leq L_1(1+|x|),
    \de
    where $\nabla b(t,x)$ stands for the derivative of $b(t,x)$ in the position $x$ and $L_1>0$ is a constant.
\end{enumerate}
\begin{enumerate}[(${\bf H}_{b}^2$)]
	\item For any $t\in[0,T]$ and $x_1, x_2\in\mR^d$
	\ce
	\|\nabla b(t,x_1)-\nabla b(t,x_2)\|\leq L_2|x_1-x_2|,
	\de
       where $L_2>0$ is a constant.
\end{enumerate}

\br
There exist a broad range of kernels which satisfy (${\bf H}_{K}^1$) (${\bf H}_{K}^2$). For example, the Riemann-Liouville kernel $K(t,s)=\frac{(t-s)_{+}^{H-\frac{1}{2}}}{\Gamma(H+\frac{1}{2})}$ for $H\in(0,1)$ satisfies (${\bf H}_{K}^1$) (${\bf H}_{K}^2$), where $\Gamma$ denotes the usual Gamma function and $(\cdot)_+:=\max\{\cdot,0\}$. Indeed, set $\a:=|\frac{1}{2}-H|$, and it holds that for any $1<\b<\frac{1}{2\a}$
\ce
\sup\limits_{t\in[0,T]}\int_0^tK(t,s)^{2\b}\dif s\leq \frac{T^{1+2\a\b}}{(1+2\a\b)\(\Gamma(H+\frac{1}{2})\)^{2\b}}, \quad H>1/2,
\de
and
\ce
\sup\limits_{t\in[0,T]}\int_0^tK(t,s)^{2\b}\dif s\leq \frac{T^{1-2\a\b}}{(1-2\a\b)\(\Gamma(H+\frac{1}{2})\)^{2\b}}, \quad H<1/2.
\de
So, $K(t,s)$ satisfies (${\bf H}_{K}^1$). For (${\bf H}_{K}^2$), note that for $s<t$
\ce
\int_0^{s} \left((s-r)^{\a}-(t-r)^{\a}\right)^2\dif r\leq T(t-s)^{2\a}, \quad H>1/2,
\de
and
\ce
&&\int_0^{s} \left((s-r)^{-\a}-(t-r)^{-\a}\right)^2\dif r=\int_0^{s}\left(\frac{(t-r)^{\a}-(s-r)^{\a}}{(s-r)^{\a}(t-r)^{\a}}\right)^2\dif r\\
&\leq&(t-s)^{2\a}\int_0^{s}\frac{1}{(s-r)^{2\a}(t-r)^{2\a}}\dif r=(t-s)^{1-2\a}\int_0^{\frac{s}{(t-s)}}\frac{1}{(v+1)^{2\a}v^{2\a}}\dif v\\
&\leq&\left\{ \begin{array}{l}C_\a(t-s)^{2\a}, ~ H\in(\frac{1}{4}, \frac{1}{2}),\\
C_\a(t-s)^{\frac{1}{2}-\a}, H\in(0,\frac{1}{4}].
\end{array}
\right.
\de
Thus, from the above deduction it follows that $K(t,s)$ satisfies (${\bf H}_{K}^2$). Besides, we remind that in the case of $0<H<1/2$, $K(t,s)$ is singular.
\er

\br
By (${\bf H}_{b,\s}^{1}$), it holds that for $t\in[0,T]$, $x,y\in\mR^{d}$
\ce
|b(t,x)-b(t,y)|\leq L_1|x-y|.
\de
\er

Under (${\bf H}_{K}^1$) (${\bf H}_{b,\s}^{1}$), by \cite[Theorem 1.1]{wang} we know that Eq.(\ref{meq1}) and Eq.(\ref{volequ}) have unique solutions $X^\e_{\cdot}, X_{\cdot}^0$, respectively. For any $\e>0$, set 
$$
Z_\cdot^\e:=\frac{X_\cdot^{\epsilon}- X_\cdot^{0}}{\sqrt{\epsilon}}, 
$$
and $Z_{\cdot}^\e$ satisfies the following stochastic Volterra equation: for any $t\in[0,T]$
\be
Z_t^\e(x)=\int_0^tK_1(t,s)\frac{b(s,X_s^{\epsilon}(x))-b(s,X_s^{0}(x))}{\sqrt{\epsilon}}\dif s+\int_0^tK_2(t,s)\sigma(s,X_s^{\epsilon}(x))\dif B_s.
\label{eq2}
\ee
Then we construct another stochastic Volterra equation:
\be
Z_{t}(x)=\int_0^tK_1(t,s)\nabla_{Z_{s}(x)}b(s,X_{s}^{0}(x))\dif s+\int_0^tK_2(t,s)\sigma(s,X_s^0(x))\dif B_s.
\label{limiequa}
\ee
And the assumptions (${\bf H}_{K}^1$) (${\bf H}_{b,\s}^{1}$) assure that Eq.(\ref{limiequa}) has a unique solution $Z_{\cdot}$. So, the CLT for Eq.(\ref{meq1}) means that as $\e$ tends to $0$,
$$
Z_{\cdot}^\e\rightarrow Z_{\cdot} ~\mbox{in distribution}.
$$

The following theorem is the main result in this paper.

\bt\label{th3}
Assume that (${\bf H}_{K}^1$) (${\bf H}_{K}^2$) (${\bf H}_{b,\s}^{1}$) and (${\bf H}_{b}^{2}$) hold. Then it holds that for any $R> 0$ and $p\geq\frac{2\b}{\b-1}$ sufficiently large,
\ce
\lim\limits_{\e\rightarrow 0}\mE\left(\sup\limits_{t\in[0,T], |x|\leq R}|Z^\e_t(x)-Z_t(x)|^p\right)=0.
\de
\et

The proof of Theorem \ref{th3} is placed in next section.

By the above theorem and the fact that the convergence in the $L^p$ sense implies the convergence in the distribution sense, we know that 
\ce
\frac{X_{\cdot}^{\epsilon}- X_{\cdot}^{0}}{\sqrt{\epsilon}}\rightarrow Z_{\cdot}
\de 
in the distribution sense, which means that $\{X_{\cdot}^{\epsilon}, \e\in(0,1)\}$ satisfies the CLT.

\section{Proof of Theorem \ref{th3}}\label{cltthproof}

In this section, we prove Theorem \ref{th3}. First of all, we make some necessary estimates. 

\bl\label{xe2p}
Under the assumptions (${\bf H}_{K}^1$) (${\bf H}_{b,\s}^{1}$), it holds that for any $p\geq\frac{2\b}{\b-1}$ and $x,y\in\mR^d$
\ce
\sup\limits_{t\in[0,T]}\mE|X^\e_{t}(x)|^{2p}\leq C(1+|x|^{2p}), \quad \sup\limits_{t\in[0,T]}\mE|X^\e_{t}(x)-X^\e_{t}(y)|^{2p}\leq C|x-y|^{2p},
\de
where the constant $C>0$ is independent of $\e$.
\el
\begin{proof}
First of all, we treat the first estimate. By the H\"older inequality, the Burkholder-Davis-Gundy inequality and (${\bf H}_{K}^1$) (${\bf H}_{b,\s}^{1}$), it holds that for any $p\geq\frac{2\b}{\b-1}$ and $0<\e<1$
\ce
\mE|X^{\epsilon}_{t}(x)|^{2p}&\leq&3^{2p-1}|x|^{2p}+3^{2p-1}\mE\left|\int_{0}^{t}K_1(t,s)b(s,X^{\epsilon}_{s}(x))\dif s\right|^{2p}\\
&&+3^{2p-1}\mE\left|\sqrt{\epsilon}\int_{0}^{t}K_2(t,s)\s(s,X^{\epsilon}_{s}(x))\dif B_s\right|^{2p}\\
&\leq&3^{2p-1}|x|^{2p}+3^{2p-1}\left(\int_{0}^{t}K_1(t,s)^{\frac{2p}{2p-1}}\dif s\right)^{2p-1}\mE\int_{0}^{t}|b(s,X^{\epsilon}_{s}(x))|^{2p}\dif s\\
&&+3^{2p-1}\mE\left|\int_{0}^{t}K_2(t,s)^2\|\s(s,X^{\epsilon}_{s}(x))\|^2\dif s\right|^{p}\\
&\leq&3^{2p-1}|x|^{2p}+3^{2p-1}\left(\int_{0}^{t}K_1(t,s)^{\frac{2p}{2p-1}}\dif s\right)^{2p-1}\mE\int_{0}^{t}|b(s,X^{\epsilon}_{s}(x))|^{2p}\dif s\\
&&+3^{2p-1}\left(\int_{0}^{t}K_2(t,s)^{\frac{2p}{p-1}}\dif s\right)^{p-1}\int_{0}^{t}\mE\|\s(s,X^{\epsilon}_{s}(x))\|^{2p}\dif s\\
&\leq&3^{2p-1}|x|^{2p}+3^{2p-1}C\int_{0}^{t}(1+\mE|X^{\epsilon}_{s}(x)|^{2p})\dif s\\
&&+3^{2p-1}C\int_{0}^{t}(1+\mE|X^{\epsilon}_{s}(x)|^{2p})\dif s,
\de
where we use the fact: 
\ce
&&\int_0^tK_1(t,s)^{\frac{2p}{2p-1}}\dif s\leq t^{1-\frac{2p}{(2p-1)\b}} \left(\int_0^tK_1(t,s)^{\b}\dif s\right)^{\frac{2p}{(2p-1)\b}}\leq C,\\
&&\int_{0}^{t}K_2(t,s)^{\frac{2p}{p-1}}\dif s\leq t^{1-\frac{p}{(p-1)\b}}\left(\int_0^tK_2(t,s)^{2\b}\dif s\right)^{\frac{p}{(p-1)\b}}\leq C.
\de
Thus, the Gronwall inequality yields that
$$
\sup\limits_{t\in[0,T]}\mE|X^{\epsilon}_{t}(x)|^{2p}\leq C(1+|x|^{2p}).
$$

Lastly, by the same deduction to that of the first estimate, we obtain the second one. 
\end{proof}

By the similar or even simpler deduction to that of Lemma \ref{xe2p}, one can obtain the following result.

\bl\label{xo2p}
For any $p\geq\frac{2\b}{\b-1}$,  it holds that for any $x,y\in\mR^d$,
\ce
\sup\limits_{t\in[0,T]}|X^0_{t}(x)|^{2p}\leq C(1+|x|^{2p}),\quad \sup\limits_{t\in[0,T]}|X^{0}_{t}(x)-X^{0}_{t}(y)|^{2p}\leq C|x-y|^{2p}.
\de
\el

\bl\label{ze2p}
Under the assumptions (${\bf H}_{K}^1$) (${\bf H}_{b,\s}^{1}$), it holds that for any $p\geq\frac{2\b}{\b-1}$
\ce
\sup\limits_{t\in[0,T]}\mE|Z^\e_{t}(x)|^{2p}\leq C(1+|x|^{2p}).
\de
\el
\begin{proof}
First of all, we observe $Z^\e$. (${\bf H}_{K}^1$), (${\bf H}_{b,\s}^{1}$), the Burkholder-Davis-Gundy inequality and the H\"older inequality imply that for any $p\geq\frac{2\b}{\b-1}$
\ce
&&\mE|Z^\e_{t}(x)|^{2p}\\
&\leq& 2^{2p-1}\mE\left|\int_0^tK_1(t,s)\frac{b(s,X_s^{\e}(x))-b(s,X_s^{0}(x))}{\sqrt{\epsilon}}\dif s\right|^{2p}\no\\
&&+2^{2p-1}\mE\left|\int_0^tK_2(t,s)\sigma(s,X_s^{\epsilon}(x))\dif B_s\right|^{2p}\no\\
&\leq& 2^{2p-1}\mE\left(\int_0^tK_1(t,s)^{\frac{2p}{2p-1}}\dif s\right)^{2p-1}\int_0^t\left|\frac{b(s,X_s^{\e}(x))-b(s,X_s^{0}(x))}{\sqrt{\epsilon}}\right|^{2p}\dif s\no\\
&&+2^{2p-1}C\mE\left|\int_0^tK_2(t,s)^2\|\sigma(s,X_s^{\epsilon}(x))\|^2\dif s\right|^{p}\no\\
&\leq& 2^{2p-1}CL_1^{2p}\int_0^t\mE|Z^{\e}_{s}(x)|^{2p}\dif s+2^{2p-1}C\left(\int_0^tK_2(t,s)^{\frac{2p}{p-1}}\dif s\right)^{p-1}\int_0^t\mE\|\sigma(s,X_s^{\epsilon}(x))\|^{2p}\dif s\no\\
&\leq& 2^{2p-1}CL_1^{2p}\int_0^t\mE|Z^{\epsilon}_{s}(x)|^{2p}\dif s+2^{2p-1}C\int_0^T(1+\mE|X_s^{\epsilon}(x)|^{2p})\dif s.
\label{zede}
\de
Thus, Lemma \ref{xe2p} and the Gronwall inequality imply the required estimate. The proof is complete.
\end{proof}

\bl\label{zepxy}
Under assumptions (${\bf H}_{K}^1$) (${\bf H}_{b,\s}^{1}$) and (${\bf H}_{b}^{2}$), it holds that for any $p\geq\frac{2\b}{\b-1}$ and $x,y\in\mR^d$
\ce
\sup\limits_{t\in[0,T]}\mE|Z^\e_{t}(x)-Z^\e_{t}(y)|^{p}\leq C(1+|x|^p)|x-y|^{p}.
\de
\el
\begin{proof}
 First of all, it holds that for any $p\geq\frac{2\b}{\b-1}$
\ce
&&\mE|Z_t^\e(x)-Z_t^\e(y)|^{p}\no\\
&\leq& 2^{p-1}\mE\left|\int_0^tK_1(t,s)\left(\frac{b(s,X_s^{\e}(x))-b(s,X_s^{0}(x))}{\sqrt{\e}}-\frac{b(s,X_s^{\e}(y))-b(s,X_s^{0}(y))}{\sqrt{\e}}\right)\dif s\right|^{p}\no\\
&&+2^{p-1}\mE\left|\int_0^tK_2(t,s)\(\sigma(s,X_s^{\epsilon}(x))-\sigma(s,X_s^{\epsilon}(y))\)\dif B_s\right|^{p}\no\\
&\leq& 2^{p-1}\mE\left(\int_0^tK_1(t,s)^{\frac{p}{p-1}}\dif s\right)^{p-1}\int_0^t\bigg|\frac{b(s,X_s^{\e}(x))-b(s,X_s^{0}(x))}{\sqrt{\e}}\\
&&\qquad\qquad -\frac{b(s,X_s^{\e}(y))-b(s,X_s^{0}(y))}{\sqrt{\e}}\bigg|^{p}\dif s\no\\
&&+2^{p-1}\mE\left|\int_0^tK_2(t,s)^2\|\sigma(s,X_s^{\epsilon}(x))-\sigma(s,X_s^{\epsilon}(y))\|^2\dif s\right|^{p/2}\no\\
&\leq& 2^{p-1}C\int_0^t\mE\bigg|\frac{b(s,X_s^{\e}(x))-b(s,X_s^{0}(x))}{\sqrt{\e}}-\frac{b(s,X_s^{\e}(y))-b(s,X_s^{0}(y))}{\sqrt{\e}}\bigg|^{p}\dif s\\
&&+2^{p-1}C\int_0^T\mE|X_s^{\epsilon}(x)-X_s^{\epsilon}(y)|^{p}\dif s.
\de

Besides, note that
\be
\frac{b(s,X_s^{\e}(x))-b(s,X_s^{0}(x))}{\sqrt{\e}}=\int_0^1\nabla_{Z_s^{\e}(x)}b(s,X_s^{0}(x)+\zeta(X_s^{\e}(x)- X_s^{0}(x)))\dif \zeta.
\label{deri}
\ee
Thus, by (${\bf H}_{b,\s}^{1}$) (${\bf H}_{b}^{2}$) it holds that
\ce
&&\left|\frac{b(s,X_s^{\e}(x))-b(s,X_s^{0}(x))}{\sqrt{\e}}-\frac{b(s,X_s^{\e}(y))-b(s,X_s^{0}(y))}{\sqrt{\e}}\right|^p\\
&\leq&2^{p-1}\bigg(\int_0^1\Big|\nabla_{Z_s^{\e}(x)}b(s,X_s^{0}(x)+\zeta(X_s^{\e}(x)- X_s^{0}(x)))\\
&&\qquad -\nabla_{Z_s^{\e}(x)}b(s,X_s^{0}(y)+\zeta(X_s^{\e}(y)- X_s^{0}(y)))\Big|\dif \zeta\bigg)^p\\
&&+2^{p-1}\bigg(\int_0^1\Big|\nabla_{Z_s^{\e}(x)}b(s,X_s^{0}(y)+\zeta(X_s^{\e}(y)- X_s^{0}(y)))\\
&&\qquad -\nabla_{Z_s^{\e}(y)}b(s,X_s^{0}(y)+\zeta(X_s^{\e}(y)- X_s^{0}(y)))\Big|\dif \zeta\bigg)^p\\
&\leq&2^{p-1}L^p_2\(2|X_s^{0}(x)-X_s^{0}(y)|+|X_s^{\e}(x)-X_s^{\e}(y)|\)^p|Z_s^{\e}(x)|^p\\
&&+2^{p-1}L^p_1|Z_s^{\e}(x)-Z_s^{\e}(y)|^p.
\de
And the H\"older inequality implies that
\ce
&&\mE\bigg|\frac{b(s,X_s^{\e}(x))-b(s,X_s^{0}(x))}{a(\e)}-\frac{b(s,X_s^{\e}(y))-b(s,X_s^{0}(y))}{a(\e)}\bigg|^{p}\\
&\leq&2^{p-1}L^p_2\left(\mE\(2|X_s^{0}(x)-X_s^{0}(y)|+|X_s^{\e}(x)-X_s^{\e}(y)|\)^{2p}\right)^{1/2}(\mE|Z_s^{\e}(x)|^{2p})^{1/2}\\
&&+2^{p-1}L^p_1\mE|Z_s^{\e}(x)-Z_s^{\e}(y)|^p\\
&\leq&2^{p-1}L^p_2C|x-y|^p(1+|x|^p)+2^{p-1}L^p_1\mE|Z_s^{\e}(x)-Z_s^{\e}(y)|^p,
\de
where in the last inequality we use Lemma \ref{xo2p}, \ref{xe2p}, \ref{ze2p}.

Finally, all the above deduction yields that  
\ce
\mE|Z_t^{\e}(x)-Z_t^{\e}(y)|^{p}\leq C(1+|x|^p)|x-y|^p+C\int_0^t\mE|Z_s^{\e}(x)-Z_s^{\e}(y)|^p\dif s.
\de
The Gronwall inequality assures the required estimate. The proof is complete.
\end{proof}

\bl\label{ttprime}
For any $p\geq 2$ sufficiently large, there exists a constant $C>0$ such that for any $t, t'\in[0,T]$ and $x\in\mR^d$
\ce
\mE|Z^\e_t(x)-Z^\e_{t'}(x)|^p\leq C(1+|x|^p)|t-t'|^{\t p},
\de
where $\t>0$ depends on $\g, \b$.
\el
\begin{proof}
For any $t<t'$, by (\ref{eq2}) we have that
\ce
Z_t^\e(x)-Z_{t'}^\e(x)&=&\int_0^t(K_1(t,s)-K_1(t',s))\frac{b(s,X_s^{\epsilon}(x))-b(s,X_s^{0}(x))}{\sqrt{\epsilon}}\dif s\\
&&+\int_t^{t'}K_1(t',s)\frac{b(s,X_s^{\epsilon}(x))-b(s,X_s^{0}(x))}{\sqrt{\epsilon}}\dif s\\
&&+\int_0^t(K_2(t,s)-K_2(t',s))\sigma(s,X_s^{\epsilon}(x))\dif B_s\\
&&+\int_t^{t'}K_2(t',s)\sigma(s,X_s^{\epsilon}(x))\dif B_s\\
&=:&I_1+I_2+I_3+I_4.
\de

For $I_1$, by (${\bf H}_{b,\s}^{1}$) (${\bf H}_{K}^{2}$) and the extended Minkowski inequality in \cite[Corollary 1.32, P.27]{ka}, it holds that
\ce
\mE|I_1|^p&\leq& L_1^p\mE\left(\int_0^t|K_1(t,s)-K_1(t',s)||Z_s^\e(x)|\dif s\right)^p\\
&\leq&L_1^p\left(\int_0^t|K_1(t,s)-K_1(t',s)|(\mE|Z_s^\e(x)|^p)^{1/p}\dif s\right)^p\\
&\leq&C(1+|x|^p)|t-t'|^{\g p}.
\de

For $I_2$, (${\bf H}_{b,\s}^{1}$) (${\bf H}_{K}^{1}$), the extended Minkowski inequality in \cite[Corollary 1.32, P.27]{ka} and the H\"older inequality imply that
\ce
\mE|I_2|^p&\leq& L_1^p\mE\left(\int_t^{t'}K_1(t',s)|Z_s^\e(x)|\dif s\right)^p\leq L_1^p\left(\int_t^{t'}K_1(t',s)(\mE|Z_s^\e(x)|^p)^{1/p}\dif s\right)^p\\
&\leq&C(1+|x|^p)\left(\int_t^{t'}K_1(t',s)\dif s\right)^p\leq C(1+|x|^p)|t-t'|^{\frac{(\b-1)p}{\b}}.
\de

By the similar deduction to that for $I_1, I_2$, we deal with $I_3, I_4$ and obtain that
\ce
\mE|I_3|^p\leq C(1+|x|^p)|t-t'|^{\g p/2}, \quad \mE|I_4|^p\leq C(1+|x|^p)|t-t'|^{\frac{(\b-1)p}{2\b}}.
\de

Finally, combining the above estimates, one can get the required result.
\end{proof}

Then by Lemma \ref{xo2p} and the similar deduction to that of Lemma \ref{ze2p}, \ref{zepxy}, \ref{ttprime}, we present the following estimate.

\bl\label{zp}
Under assumptions (${\bf H}_{K}^1$) (${\bf H}_{K}^2$) (${\bf H}_{b,\s}^{1}$) and (${\bf H}_{b}^{2}$), it holds that 
\ce
&&\sup\limits_{t\in[0,T]}\mE|Z_{t}(x)|^{p}\leq C(1+|x|^p), \quad x\in\mR^d, p\geq\frac{2\b}{\b-1},\\
&&\sup\limits_{t\in[0,T]}\mE|Z_{t}(x)-Z_{t}(y)|^{p}\leq C(1+|x|^p)|x-y|^{p}, \quad x,y\in\mR^d, p\geq\frac{2\b}{\b-1},\\
&&\mE|Z_t(x)-Z_{t'}(x)|^p\leq C(1+|x|^p)|t-t'|^{\t p}, \quad t, t'\in[0,T], x\in\mR^d, p\geq 2.
\de
\el

The following lemma is important to prove the central limit theorem for Eq.(\ref{meq1}).

\bl\label{cl}
Under assumptions (${\bf H}_{K}$) (${\bf H}_{b,\s}^{1}$) and (${\bf H}_{b}^{2}$), it holds that for any $p\geq\frac{2\b}{\b-1}$
\ce
\sup\limits_{t\in[0,T]}\mE\left|Z^\e_{t}(x)-Z_{t}(x)\right|^{p}\leq C(1+|x|^{2p})\epsilon^{p/2}, \quad x\in\mR^d,
\de
where the constant $C>0$ is independent of $\e$.
\el
\begin{proof}
We begin with (\ref{eq2}) (\ref{limiequa}). For any $p\geq\frac{2\b}{\b-1}$, it holds that 
\ce
&&\mE|Z_t^\e(x)-Z_{t}(x)|^p\\
&\leq& 2^{p-1}\mE\left|\int_0^tK_1(t,s)\left(\frac{b(s,X_s^{\epsilon}(x))-b(s,X_s^{0}(x))}{\sqrt{\epsilon}}-\nabla_{Z_{s}(x)}b(s,X_{s}^{0}(x))\right)\dif s\right|^p\\
&&+2^{p-1}\mE\left|\int_0^tK_2(t,s)(\sigma(s,X_s^{\epsilon}(x))-\sigma(s,X_s^0(x)))\dif B_s\right|^p\\
&\leq& 2^{p-1}\mE\left(\int_0^tK_1(t,s)^{\frac{p}{p-1}}\dif s\right)^{p-1}\int_0^t\bigg|\frac{b(s,X_s^{\epsilon}(x))-b(s,X_s^{0}(x))}{\sqrt{\epsilon}}\\
&&\qquad\qquad -\nabla_{Z_{s}(x)}b(s,X_{s}^{0}(x))\bigg|^p\dif s\\
&&+2^{p-1}\mE\left|\int_0^tK_2(t,s)^2\|\sigma(s,X_s^{\epsilon}(x))-\sigma(s,X_s^0(x))\|^2\dif s\right|^{p/2}\\
&\leq& 2^{p-1}C\mE\int_0^t\left|\frac{b(s,X_s^{\epsilon}(x))-b(s,X_s^{0}(x))}{\sqrt{\epsilon}}-\nabla_{Z_{s}(x)}b(s,X_{s}^{0}(x))\right|^p\dif s\\
&&+2^{p-1}C\mE\int_0^t\|\sigma(s,X_s^{\epsilon}(x))-\sigma(s,X_s^0(x))\|^p\dif s\\
&=:&I_1+I_2.
\de
For $I_1$, (\ref{deri}) (${\bf H}_{b,\s}^{1}$) and (${\bf H}_{b}^{2}$) imply that
\ce
I_1&\leq&4^{p-1}C\mE\int_0^t\left|\int_0^1\nabla_{Z_s^\e(x)}b(s,X_s^{0}(x)+\zeta(X_s^{\epsilon}(x)- X_s^{0}(x)))\dif \zeta-\nabla_{Z^\e_{s}(x)}b(s,X_{s}^{0}(x))\right|^p\dif s\\
&&+4^{p-1}C\mE\int_0^t\left|\nabla_{Z^\e_{s}(x)}b(s,X_{s}^{0}(x))-\nabla_{Z_{s}(x)}b(s,X_{s}^{0}(x))\right|^p\dif s\\
&\leq&4^{p-1}CL_2^p\mE\int_0^t\left|X_s^{\epsilon}(x)- X_s^{0}(x)\right|^p|Z^\e_{s}(x)|^p\dif s+4^{p-1}CL_1^p\mE\int_0^t\left|Z^\e_{s}(x)-Z_{s}(x)\right|^p\dif s\\
&\leq&4^{p-1}CL_2^p\e^{p/2}\int_0^t\mE|Z^\e_{s}(x)|^{2p}\dif s+4^{p-1}CL_1^p\int_0^t\mE\left|Z^\e_{s}(x)-Z_{s}(x)\right|^p\dif s.
\de
We deal with $I_2$. By (${\bf H}_{b,\s}^{1}$), it holds that
\ce
I_2\leq 2^{p-1}CL_1^p\mE\int_0^t|X_s^{\epsilon}(x)-X_s^0(x)|^p\dif s\leq 2^{p-1}CL_1^p\e^{p/2}\int_0^t\mE|Z^\e_{s}(x)|^{p}\dif s.
\de
Combining all the above deduction, we conclude that
\ce
&&\mE|Z_t^\e(x)-Z_{t}(x)|^p\no\\
&\leq& 4^{p-1}CL_2^p\e^{p/2}\int_0^t\mE|Z^\e_{s}(x)|^{2p}\dif s+4^{p-1}CL_1^p\int_0^t\mE\left|Z^\e_{s}(x)-Z_{s}(x)\right|^p\dif s\no\\
&&+2^{p-1}CL_1^p\e^{p/2}\int_0^t\mE|Z^\e_{s}(x)|^{p}\dif s.
\label{zez}
\de

Lastly, by Lemma \ref{ze2p} and the Gronwall inequality one can get that 
\ce
\sup\limits_{t\in[0,T]}\mE|Z_t^\e(x)-Z_{t}(x)|^p\leq C(1+|x|^{2p})\e^{p/2}.
\de
The proof is complete.
\end{proof}

{\bf Proof of Theorem \ref{th3}.} First of all, we take any sequence $\{\e_n, n\in\mN\}$ satisfying $\e_n\downarrow0$ as $n\rightarrow\infty$ and construct the following process:
\ce
{\bf Z}(r,t,x)=\left\{\begin{array}{l}
Z_t(x), \quad r=0 \\
Z_t^{\e_n}(x)+\e^{-1}_{n+1}\left(r-\e_nT\right)\left[Z_t^{\e_{n+1}}(x)-Z_t^{\e_{n}}(x)\right], \quad \e_{n+1}T<r \leqslant \e_nT,\quad n \in \mathbb{N}.
\end{array}\right.
\de
For any $R>0$, Lemma \ref{zepxy}, \ref{ttprime}, \ref{zp} imply that for any $p\geq\frac{2\b}{\b-1}$, $t, t'\in[0,T]$, and $x,y \in D_R:=\{x\in\mR^d; |x|\leq R\}$ 
\ce
\mE|{\bf Z}(r,t,x)-{\bf Z}(r',t',y)|^p\leq C(|r-r'|^{\eta p}+|t-t'|^{\eta p}+|x-y|^{\eta p}),
\de
where the constant $\eta>0$ depends on $\b,\g,\t$. Thus, for $p$ sufficiently large, by Kolmogorov's continuity criterium, there is a $p$-order integrable random variable $\xi$ such that
$$
\sup _{t\in[0,T], |x| \leqslant R}\left|{\bf Z}(r, t, x)-{\bf Z}\left(r^{\prime}, t, x\right)\right| \leqslant \xi\left|r-r^{\prime}\right|^\lambda, \quad \text { a.s. }
$$
where $\lambda \in\left(0, \eta-\frac{d+2}{p}\right)$. Especially, we take $r^{\prime}=0, r=\e_nT$, and obtain that
$$
\mathbb{E}\left(\sup _{t \in[0,T], |x| \leqslant R}\left|Z_t^{\e_n}(x)-Z_t(x)\right|^p\right) \leqslant T^{\lambda p}\mathbb{E}|\xi|^p\e_n^{\lambda p},
$$
which yields the desired convergence.

\section{Application}

In this section, we apply our result to stochastic Volterra equations with the kernels of fractional Brownian motions with the Hurst parameter $H\in(0, 1)$.

First of all, we recall some basics about fractional Brownian motions (c.f. \cite{cd, du}). Fix $T>0$. Let $\{B_t^H, t\in[0,T]\}$ be a fractional Brownian motion with the Hurst index $H\in(0,1)$, which is a centered Gaussian process with the following covariance 
$$
R_H(s,t)=\frac{V_H}{2}(s^{2H}+t^{2H}-|t-s|^{2H}),
$$
where $V_H=\frac{\Gamma(2-2H)\cos(\pi H)}{\pi H(1-2H)}$, and $\Gamma$ denotes the usual Gamma function. It is known that the fractional Brownian motion $B^H$ has the representation in law:
$$
B_t^H=\int_0^t K_H(t,s)\dif B_s,
$$
where $K_H(t,s)$ is the square root of the covariance operator, that is
$$
R_H(s,t)=\int_0^T K_H(s,r)K_H(t,r)\dif r.
$$
More precisely,
$$
K_H(t,r)=\frac{(t-r)^{H-\frac{1}{2}}}{\Gamma(H+\frac{1}{2})}F(\frac{1}{2}-H, H-\frac{1}{2}, H+\frac{1}{2}, 1-\frac{t}{r})I_{[0,t)}(r),
$$
where $F$ is the Gauss hypergeometric function.

Next, we justify that $K_H(t,s)$ satisfies (${\bf H}_{K}^1$) (${\bf H}_{K}^2$). Note that $0\leq K_H(t,s)\leq Cs^{-|H-1/2|}(t-s)^{-(1/2-H)_+}I_{s<t}$. Set $H_0:=|H-1/2|$ and it holds that for $1<\b<\frac{1}{2H_0}$,
\ce
\sup\limits_{t\in[0,T]}\int_0^tK_H(t,s)^{2\b}\dif s\leq CT^{1-2\b H_0}{\bf B}(1-2\b H_0,1-2\b H_0), \quad H<1/2,
\de
where ${\bf B}$ is the usual Beta function, and 
$$
\sup\limits_{t\in[0,T]}\int_0^tK_H(t,s)^{2\b}\dif s\leq C\frac{T^{1-2\b H_0}}{1-2\b H_0}, \quad H>1/2.
$$
So, $K_H(t,s)$ satisfies (${\bf H}_{K}^1$). Besides, for $t<t'$
\ce
\int_0^{t\land t'}|K_H(t',s)-K_H(t,s)|^2\dif s&\leq&\int_0^TK_H(t',s)K_H(t',s)\dif s-2\int_0^TK_H(t',s)K_H(t,s)\dif s\\
&&+\int_0^TK_H(t,s)K_H(t,s)\dif s\\
&=&R(t',t')-2R(t',t)+R(t,t)\leq C|t-t'|^{2H}.
\de
That is, $K_H(t,s)$ satisfies (${\bf H}_{K}^2$).

Finally, we take $K_1(t,s)=K_2(t,s)=K_H(t,s)$ and consider the following stochastic Volterra equation: for $0<\e<1$
\ce
X^\e_t(x)=x+\int_0^tK_H(t,s)b(s,X^\e_s(x))\dif s+\sqrt{\e}\int_0^tK_H(t,s)\s(s,X^\e_s(x))\dif B_s, \quad t\in[0,T].
\de
Assume that $b,\s$ satisfy (${\bf H}_{b,\s}^{1}$) and (${\bf H}_{b}^{2}$). Then the above equation has a unique solution $X^\e$. Moreover, by Theorem \ref{th3}, we obtain that 
\ce
\frac{X^\e-X^0}{\sqrt{\e}}\overset{L^p}{\rightarrow} Z,
\de
where $X^0$ is the solution of the following Volterra equation
\ce
X^0_t(x)=x+\int_0^tK_H(t,s)b(s,X^0_s(x))\dif s,
\de
and $Z$ is the solution of the following stochastic Volterra equation
\ce
Z_{t}(x)=\int_0^tK_H(t,s)\nabla_{Z_{s}(x)}b(s,X_{s}^{0}(x))\dif s+\int_0^tK_H(t,s)\sigma(s,X_s^0(x))\dif B_s.
\de

\end{document}